\documentclass[11pt]{article}
\usepackage{color}
\usepackage[affil-it]{authblk}
\usepackage{amssymb}
\usepackage{amsmath,mathtools}
\usepackage{amsfonts}
\usepackage{amsthm}
\usepackage{dsfont}
\usepackage{amssymb}
\usepackage{natbib}
\usepackage{mathpazo}
\usepackage{fullpage}
\usepackage{enumerate}
\usepackage[colorlinks=true,linkcolor=blue,citecolor=blue,urlcolor=blue]{hyperref}

\usepackage{titlesec}
\titleformat{\section}[hang]{\Large\bfseries\filright}{\thesection.}{.5em}{}
\titleformat{\subsection}[hang]{\large\bfseries\filright}{}{0em}{}
\titleformat{\subsubsection}[block]{\bfseries}{}{0em}{}

\newtheorem{definition}{Definition}
\newtheorem{theorem}{Theorem}

\newtheorem{lemma}[theorem]{Lemma}
\newtheorem{remark}{Remark}

\newcommand{\M}{\mathcal{M}}
\newcommand{\N}{\mathcal{N}}
\newcommand{\Hi}{\mathcal{H}}

\newcommand{\Tr}{\operatorname{Tr}}

\title{Jensen's inequality for partial traces in von Neumann algebras}
\author [1, 2]{Mizanur Rahaman}
\author [2]{Lyudmila Turowska}
\affil [1]{Wallenberg Centre for Quantum Technology, Chalmers University of Technology}
\affil [2]{Department of Mathematical Sciences, Chalmers University of Technology and University of Gothenburg}
\date{\today}

\begin{document}
\maketitle

\begin{abstract}Motivated by a recent result on finite-dimensional Hilbert spaces, we prove a Jensen's inequality for partial traces in semifinite von Neumann algebras. We also prove a similar inequality in the framework of general (non-tracial) von Neumann algebras.
\end{abstract}
\section{Introduction}
Jensen’s inequality is a very useful inequality in convex analysis and probability
theory. It provides a relationship between convex functions and expectations (or integrals). Extending this connection, Davis \cite{davis}, and later Choi \cite{choi1974} proved Jensen's inequality for positive maps and operator convex functions. Also, Hansen–Pedersen showed an operator inequality (\cite{Hansen-Pedersen}) which is a powerful and elegant generalization of Jensen’s inequality to the setting of operator algebras, particularly involving operator convex functions. These results play a crucial role in functional analysis and have wide applications in probability and  quantum information theory.

In a recent article (\cite{Carlen-Frank=Larson}), Carlen, Frank, and Larson established a Jensen's inequality for partial traces on finite-dimensional Hilbert spaces. More precisely, they showed the following theorem.
\begin{theorem}[Carlen-Frank-Larson, 2025]\label{CFL theorem}
For two finite-dimensional Hilbert spaces $\mathcal{H}_1, \mathcal{H}_2$ and a self-adjoint operator $H$ acting on $\mathcal{H}_1\otimes \mathcal{H}_2$ whose spectrum lies within an interval $I$, we have 
\[\Tr_2f(\Tr_1(\rho\otimes \mathbb{1})^{\frac{1}{2}}H(\rho\otimes \mathbb{1})^{\frac{1}{2}})\leq \Tr_1 (\rho^\frac{1}{2} \Tr_2 f(H)\rho^\frac{1}{2}),\]
for every convex function $f$ defined on $I$ and every density matrix (positive semi-definite matrix with trace $1$) $\rho$ on $\mathcal{H}_1$. 
\end{theorem}
This inequality helped to establish some results about eigenvalue asymptotics of operators with homogeneous potentials. 


Owing to the general interest in Jensen's inequality the aim of this article is to extend this new result in the framework of von Neumann algebras. Some notable works have already been done in generalizing tracial version of Jensen's inequality in operator algebras: pioneered by Brown and Kosaki \cite{brown-kosaki} who proved a trace Jensen's inequality in semifinite von Neumann algebras. Later Petz showed similar inequality involving contractive positive maps \cite{Petz}. Harada-Kosaki \cite{harada-kosaki} further improved both the previous results involving semi-bounded self-adjoint operators. We also note that there are many variations of Jensen's inequality in different contexts, and for the reader's convenience, we give some references here \cite{massey, farenick2007, hansen-pedersen2}. It is noteworthy that such an inequality for partial traces on operator algebraic contexts has not been established. 

We state our main theorem below.
\begin{theorem} [Jensen's inequality for partial traces in von Neumann algebras]\label{main theorem}
    Let $(\M_1, \tau_1)$ and $(\M_2, \tau_2)$ be two von Neumann algebras with the corresponding normal semifinite traces. Let $H\in \M_1\bar{\otimes} \M_2$ be a self-adjoint element, and $f:\mathbb{R}\rightarrow \mathbb{R}$ a continuous convex function. For $a\in L^2(\M_1,\tau_1)$, with $\tau_1(a^*a)=1$, we have 

\begin{equation}\label{eq-main}
    \tau_2 f (\tau_1\otimes id)[(a^*\otimes 1)H (a\otimes 1)]\leq \tau_1 [a^*( id\otimes \tau_2) f(H)a], 
    \end{equation}
     whenever both sides are defined.
     
\noindent Moreover, if $\tau_1(a^*a)\leq 1$, and $f(0)=0$, then the stated inequality \ref{eq-main} still holds whenever both sides are defined.
\end{theorem}

See remark \ref{remark} for what we mean by both sides of the inequality being defined. 

We remark that when $\M_1=B(\mathcal{H}_1)$ and $
\M_2=B(\mathcal{H}_2)$, where $\mathcal{H}_1,$ and $\mathcal{H}_2$ are finite-dimensional Hilbert spaces, the inequality above is slightly stronger than the statement in Theorem \ref{CFL theorem}, as in their work, the authors need the square root of a density operator in the conjugation. For the finite-dimensional case, any operator is trace-class and hence the extra assumption on the element $a$ in our theorem is not needed. We would like to further point out that we proved a similar inequality in Theorem \ref{main theorm states} for general (non-tracial) von Neumann algebras that involve normal states but here we require an operator convex function for the inequality.

\section{Preliminary on semifinite von Neumann algebras}
We refer to \cite[Chapter 1]{pisier-xu} and \cite[Chapter V.2]{takesaki} for some preliminary contexts on von Neumann algebras and we give a brief outline of semifinite traces below. 
Let  $\mathcal M$ be a von Neumann algebra acting on $\Hi$, and write $\mathcal M_+$ for the cone of positive elements of $\mathcal M$. If $\mathcal M$ is equipped with a normal faithful semifinite trace, we call $\mathcal M$ tracial.  Recall that  a trace  on $\mathcal M$ is $\tau:\mathcal M_+\to[0,\infty]$ satisfying
\begin{enumerate}
\item $\tau(x+y)=\tau(x)+\tau(y)$, $x$, $y\in \mathcal M_+$;
\item $\tau(\lambda x)=\lambda\tau(x)$, $\lambda\geq 0$, $x\in \mathcal M_+$;
\item $\tau(xx^*)=\tau(x^*x)$, $x\in\mathcal M$.
\end{enumerate}
$\tau$ is said to be faithful if $\tau(x)=0$ implies $x=0$, and $\tau$ is semifinite if for any $x\in\mathcal M_+$ there exists a non-zero $y\in\mathcal M_+$ such that $y\leq x$ and $\tau(y)<\infty$.  

Given any tracial von Neumann algebra, let $S_+$ be the set of all $x\in\mathcal M_+$ such that $\tau({\rm supp} x)<\infty$, where ${\rm supp} x$ denotes the support of $x$, that is the least projection $p$ such that $xp=x$.  Let $S$ be the linear span of $S_+$. Then for any $x\in S$ one has $\tau(|x|)<\infty$. The completion of $S$ with respect $\|x\|_1:=\tau(|x|)$ is the non-commutative $L^1$-space $L^1(\mathcal M,\tau)$ 
 associated with $(\mathcal M,\tau)$. The space $\mathcal M\cap L^1(\mathcal M,\tau)$ is dense in $L^1(\mathcal M,\tau)$, the trace $\tau$ is finite on $\mathcal M_+\cap L^1(\mathcal M,\tau)$ and extends to a functional on $L^1(\mathcal M,\tau)$.  
 
 The elements of $L^1(\mathcal M,\tau)$ can be considered as closed densely defined operators on $\mathcal{H}$ affiliated with $\mathcal M$. An affiliated element is called measurable with respect to $\tau$ if $\tau(e_\lambda(|x|))<\infty$ for some $\lambda>0$, $e_\lambda$ is the characteristic function of the set $(\lambda, \infty)$.
 Let $L^0(\mathcal M,\tau)$ denote the space of measurable elements. It is known that $L^0(\mathcal M,\tau)$ is a $*$-algebra and  $\tau$ can be extended to the positive part of $L^0(\mathcal M,\tau)$. Then $L^1(\mathcal M,\tau)=\{x\in L^0(\mathcal M,\tau): \tau(|x|)<\infty\}$ is a Banach space with the norm $\|x\|_1=\tau(|x|)$.  
For  $x\in L^0(\M,\tau)$ recall the generalized singular numbers 
$$\mu_t(x)={\rm inf}\{\lambda>0: \tau(e_\lambda(|x|))\leq t\}, \quad t>0.$$
It is known that $\tau$ can be extended to $L^0(\M,\tau)_+$, by letting 
$\tau(x)=\int_0^\infty\mu_t(x)dt$, $x\in L^0(\M,\tau)_+$.
Similar to $L^1(\M,\tau)$ one defines the non-commutative $L^2$-space, $L^2(\M,\tau)$, as all elements $x\in L^0(\M,\tau)$ such that $\tau(x^*x)<\infty$. One has that for any $x$, $y\in L^2(\M,\tau)$, $xy^*\in L^1(\M,\tau)$ and $\langle x,y\rangle=\tau(xy^*)$ is an inner product making $L^2(\M,\tau)$ a Hilbert space. Writing $\|x\|_2$ for the corresponding norm we have $\|xa\|_2^2=\tau(a^*x^*xa)\leq \|x\|^2\|a\|_2^2$ for any $x\in \M$, $a\in L^2(\M,\tau)$. 


\section{Trace Jensen's inequality for positive linear maps}
Recall that a function $f$ defined on an interval $I$ is convex if $f(\lambda x+ (1-\lambda) y)\leq \lambda f(x)+(1-\lambda) f(y)$ whenever $\lambda\in [0,1]$ and $x,y\in I$. 

We start this section by stating Petz's result (see Theorem A and Theorem B in \cite{Petz}).
\begin{theorem}[Petz. 1987]\label{Petz thm}
    Let $\M, \mathcal{N}$ be von Neumann algebras and $f:\mathbb{R}^+\rightarrow \mathbb{R}$ be a continuous convex function. Assume that $\tau$ is a normal semifinite trace on $\M$ and $\Phi:\mathcal{N}\rightarrow \M$ is a unital positive map. Then for any $a\in \mathcal{N}_+$ one has
    \[ \tau(f(\Phi(a)))\leq \tau( \Phi(f(a))),\]
    provided both sides exist.
    
    Furthermore, if we further assume that $f(0)=0$, then the above inequality remains valid for contractive  positive $\Phi$.
\end{theorem}

Note that a positive unital linear mapping on a C$^*$-algebra is necessarily a contraction (\cite{paulsen-book}). In this section we show that the assertion of Petz's result holds even when we replace the positive element in the statement of Theorem \ref{Petz thm} by a self-adjoint element. This is a crucial step and the main technical ingredient in proving our main result given in the next section. 

Following \cite{brown-kosaki} we make the following definition
\begin{definition} \label{def-Jordan}Let $x$ be a selfadjont operator in $\M$ with the Jordan decomposition $x=x_+-x_-$. We say that $\tau(x)$ is defined if $\tau(x_+)<+\infty$ or $\tau(x_-)<+\infty$, in this case we set $\tau(x)=\tau(x_+)-\tau(x_-)\in[-\infty, +\infty]$.
\end{definition}

We have 
\begin{theorem}\label{improved Petz thm}
    Let $\M, \mathcal{N}$ be von Neumann algebras and $f:\mathbb{R}\rightarrow \mathbb{R}$ be a continuous convex function. Assume that $\tau$ is a normal semifinite trace on $\M$ and $\Phi:\mathcal{N}\rightarrow \M$ is a unital positive map. Then for any self-adjoint $x\in \N$ we have
    \[ \tau(f(\Phi(x)))\leq \tau (\Phi(f(x))),\]
    provided both sides exist.

    Furthermore, the above inequality remains valid if we have a contractive positive map $\Phi$ and a continuous convex function with $f(0)=0$. 
\end{theorem}
Note that the above statement not only improves Petz's result, but also gives a generalization of one of the main results of \cite{harada-kosaki}, as in Theorem 7 in \cite{harada-kosaki} they prove a trace Jensen's inequality for maps of the form $x\rightarrow a^*xa$, where $a$ is a contraction.

We need a few lemmas, which are generalizations of Lemma 4 and Lemma 5 given in \cite{harada-kosaki}. Indeed, in \cite{harada-kosaki} the following lemma was stated for the special map $\Phi(x)=a^*xa$, for a contraction $a\in \M$. Here we show that the assertion holds for a 
positive contractive map. In the following we will use the concept of spectral pre-order. 
\begin{definition}
    Let $a, b\in \M$ be two self-adjoint elements. Then we say that $a\lesssim b$ if $e_{(s, \infty)}(a)$, is 
    equivalent, in the Murray-von Neumann sense, to a subprojection of $e_{(s, \infty)}(b)$ for every real number $s$, where $e_I(z)$ is the spectral projection of any self-adjoint element $z\in \M$ corresponding to the subset $I\subseteq \mathbb{R}$. 
\end{definition}
We remark that if $a\lesssim b$ and $\tau$ is a normal semifinite trace on $\M$ then $\tau(a)\leq \tau(b)$ (see \cite[Lemma 3]{brown-kosaki}.

The following Lemma was proved in \cite{brown-kosaki} for maps $\Phi: B(\Hi)\to B(\Hi)$, $\Phi(x)=a^*xa$, where $a\in B(\Hi)$ is a contraction. It is easily generalized to any contractive positive normal map on a von Neumann algebra $\M$. 
\begin{lemma}\label{contr}
Let $\Phi: \M\rightarrow \M $ be a unital positive map and suppose $\M\subseteq B(\mathcal{H})$. Then for any continuous convex function  $f:{\mathbb R}\to{\mathbb R}$ and  a self-adjoint $x\in \M$ we have
\[f(\langle\Phi(x)\xi, \xi\rangle)\leq \langle \Phi(f(x))\xi, \xi \rangle,\]
for any unit vector $\xi\in \mathcal{H}$. 

Furthermore, if $\Phi$ is assumed to be  a positive contraction, then the above inequality holds for any continuous convex function $f:\mathbb R\to\mathbb R$ with the assumption $f(0)=0$.

\end{lemma}
\begin{proof}
 Note that a positive unital map on $\M$ is necessarily contractive. Let  $\Phi:\M\rightarrow \M $ be a positive contraction.
As $g\mapsto \langle \Phi(g(x))\xi,\xi\rangle$ is a positive bounded linear functional on $C_0(\mathbb R)$, by Riesz theorem there exists a Borel measure $\mu_\Phi$ such that 
$$\langle \Phi(g(x))\xi,\xi\rangle=\int_{\mathbb R} g(t)d\mu_\Phi(t).$$
As the spectrum $\sigma(x)$ is a compact subset of $\mathbb R$ and $g(x)=0$ if $\text{supp }g\cap\sigma(x)=\emptyset$, $\text{supp }\mu_\Phi$ is compact.   Taking $g\in C_0(\mathbb R)$, $g=1$ on a compact containing $\sigma(x)$ and $\text{supp }\mu_\Phi$, and using the fact that $\Phi$ is a positive contraction,  we obtain $\mu_\Phi(\mathbb  R)=\langle\Phi(1)\xi,\xi\rangle\leq 1$. Let $\nu_\Phi=\mu_\Phi+(1-\langle\Phi(1)\xi,\xi\rangle)\delta_0$, where $\delta_0$ is the Dirac measure supported in $\{0\}$. Then $\nu_\Phi$ is a probability measure. Now suppose that  $f(0)=0$. As the functions $f(t)$ and $g(t)=t$ vanish in $0$ we obtain from the classical Jensen's inequality that
\begin{eqnarray*}
&&\langle \Phi(f(x))\xi,\xi\rangle=\int_{\mathbb R} f(t)d\mu_\Phi(t)=\int_{\mathbb R} f(t)d\nu_\Phi(t)\\&&\geq f\left(\int_{\mathbb R} td\nu_\Phi(t)\right)=f\left(\int_{\mathbb R} td\mu_\Phi(t)\right)=f(\langle \Phi(x)\xi,\xi\rangle).
\end{eqnarray*}

Note that if $\Phi$ is unital, then $\mu_\Phi$ is a probability measure and the rest of the proof follows for any continuous convex function $f$ without the assumption $f(0)=0$. 
\end{proof}

\begin{lemma}\label{key lemma}
Let $\Phi: \M\rightarrow \M $ be a unital  positive map. For any continuous convex function $f$ and  a self-adjoint $x\in \M$, let $I\subset \mathbb{R}$ be an interval on which $f(t)$ is monotone, that is, either increasing or decreasing. And set $p=e_I(\Phi(x))$ to be the spectral projection of $\Phi(x)$ on the interval $I$.
    \begin{enumerate}[(a)]

            \item If $f(t)\geq 0$ on the interval $I$, then we have $p \Phi(f(x))p\geq 0$ and moreover $$p \Phi(f(x))p\gtrsim pf(\Phi(x))p.$$
            \item If $f(t)\leq 0$ on the interval $I$, then the negative part $(p\Phi(f(x))p)\_$ of the Jordan decomposition satisfies  $$-p f(\Phi(x))p \gtrsim(p\Phi(f(x))p)\_.$$
    
    \end{enumerate}
    Furthermore, if $\Phi$ is assumed to be a positive contraction, then the assertions of $(a)$ and $(b)$ follow for any continuous convex functions $f$ with $f(0)=0$.
\end{lemma}
 \begin{proof}

 (a)
                  To see how to prove  $p \Phi(f(x))p\geq 0$,  we note that if we choose $\xi\in \text{Range}(p)$, then using Lemma \ref{contr} yields $\langle\Phi(f(x)) \xi, \xi\rangle \geq f(\langle\Phi(x)\xi, \xi\rangle)\geq 0$. The last inequality follows because $f$ was assumed to be positive on $I$.


We have to prove that for any $s>0$ there exists a projection $q_s$ such that
         $$e_{(s,+\infty)}(pf(\Phi(x))p)\sim q_s\leq e_{(s,+\infty)}(p\Phi(f(x))p).$$
         For $s>0$ we write $[f>s]\cap I$ for $\{t\in I: f(t)>s\}$.
         We observe first that $e_{(s,+\infty)}(pf(\Phi(x))p)=e_{[f>s]\cap I}(\Phi(x))$, as $\chi_I(t)f(t)>s$ if and only if $t\in I$ and $f(t)>s$. 
         Assuming that $e_{[f>s]\cap I}(\Phi(x))$  is non-zero, in particular, $[f>s]\cap I\ne\emptyset$, and taking 
         $\xi$  a unit vector in $\text{Range}(e_{[f>s]\cap I}(\Phi(x))$, we obtain using monotonicity of $f$ on $I$, that 
         $\langle \Phi(x)\xi,\xi\rangle\in \overline{[f> s]\cap I}$. In fact, letting  $q=e_{[f>s]\cap I}(\Phi(x))$, the spectrum $\sigma(q\Phi(x)q)$ is in the closure of $[f>s]\cap I$ which is a closed interval, say $[\alpha,\beta]\subset\overline{I}$,  giving  $\alpha\leq \langle\Phi(x)\xi,\xi\rangle\leq\beta$. Assume that $f$ is increasing on $I$ (the case when $f$ is decreasing is treated similarly). Then either $f(\alpha)=s$ or $f(\alpha)>s$; in the latter case  $f(\langle\Phi(x)\xi,\xi\rangle)\geq f(\alpha)>s$ and hence $$s<f(\langle\Phi(x)\xi,\xi\rangle)=f(\langle p\Phi(x)p\xi,\xi\rangle);$$
         and in the former case if 
         $f(\langle\Phi(x)\xi,\xi\rangle)=f(\alpha)=s$, then $\langle q\Phi(x)q\xi,\xi\rangle=\langle\Phi(x)\xi,\xi\rangle\leq\alpha$, as $f(t)>s$ for all $t>\alpha$, $t\in I$; on the other hand, $q\Phi(x)q\geq\alpha$, giving $\Phi(x)\xi=q\Phi(x)q\xi=\alpha\xi$ and $\xi=e_{[f=s]}(\Phi(x))\xi$, a contradiction. Therefore, $$f(\langle p\Phi(x)p\xi,\xi\rangle)=f(\langle\Phi(x)\xi,\xi\rangle)>f(\alpha)=s.$$
         Applying Lemma \ref{contr}, we get $$\langle p\Phi(f(x))p\xi,\xi\rangle>s$$ giving  $e_{(-\infty, s]}(p\Phi(f(x))p)\xi=0$ and hence $$e_{(s,+\infty)}(pf(\Phi(x))p)\wedge e_{(-\infty,s]}(p\Phi(f(x))p)=0.$$
         If $e_{(s,+\infty)}(pf(\Phi(x))p)=e_{[f>s]\cap I}(\Phi(x))=0$, the latter holds trivially.
         
Now using Kaplanski's formula $p_1\vee p_2-p_1\sim p_2-p_1\wedge p_2$ for any pair of projections $(p_1,p_2)$,
we obtain 
         \begin{eqnarray*}
         &&e_{(s,+\infty)}(pf(\Phi(x))p) =e_{(s,+\infty)}(pf(\Phi(x))p)-e_{(s,+\infty)}(pf(\Phi(x))p)\wedge e_{(-\infty,s]}(p\Phi(f(x))p)\\&&\sim e_{(s,+\infty)}(pf(\Phi(x))p)\vee e_{(-\infty, s]}(p(\Phi(f(x)))p)- e_{(-\infty,s]}(p\Phi(f(x))p)\\&&\leq I- e_{(-\infty,s]}(p\Phi(f(x))p)=e_{(s,+\infty)}(p\Phi(f(x))p).
         \end{eqnarray*}

   (b) is obtained using similar arguments.      
          
 \end{proof}

\textbf{Proof of Theorem \ref{improved Petz thm} :}
We follow the same line of argument as given in \cite{harada-kosaki}. Without loss of generality we consider $f(t)$ is decreasing at the origin. We assume that $f(t)$ is
\begin{enumerate}
    \item positive and decreasing on $I_1=(-\infty, 0),$
    \item negative and decreasing on $I_2=[0, t_1),$
    \item negative and increasing on $I_3=[t_1, t_2),$
    \item positive and increasing on $I_4=[t_2, \infty)$. 
\end{enumerate}
Some of these intervals may be empty, for example, if $f(t)\geq 0$, then $I_2=I_3=\varnothing$. Now for the self-adjoint element $x\in \M$ we set $p_i=e_{I_i}(\Phi(x))$ and let $E(y)=\sum_{i=1}^4 p_i y p_i$. 
We can write $f(\Phi(x))=E(f(\Phi(x)))=\sum_{i}p_i f(\Phi(x)) p_i$ in the matrix form:

\[f(\Phi(x))=\begin{bmatrix}
p_1 f(\Phi(x)) p_1 & 0 & 0 & 0\\
0 & p_2 f(\Phi(x)) p_2 & 0 & 0\\
0 & 0 & p_3 f(\Phi(x)) p_3 & 0\\
0 & 0 & 0 & p_4 f(\Phi(x)) p_4
\end{bmatrix}.\]
Note that for $i=1, 4$, $p_i f(\Phi(x)) p_i\geq 0$ and for $i=2, 3$, we have $p_i f(\Phi(x)) p_i\leq 0$.

However, $\Phi(f(x))$ is not diagonal, and we have 

\[E(\Phi(f(x)))=\begin{bmatrix}
p_1 \Phi(f(x)) p_1 & 0 & 0 & 0\\
0 & p_2 \Phi(f(x)) p_2 & 0 & 0\\
0 & 0 & p_3 \Phi(f(x)) p_3 & 0\\
0 & 0 & 0 & p_4 \Phi(f(x)) p_4
\end{bmatrix}.\]

Note that using part $(a)$ of the Lemma \ref{key lemma} we have $p_1 \Phi(f(x)) p_1\geq 0$ and $p_4 \Phi(f(x)) p_4\geq 0$.

 The Jordan decomposition of $f(\Phi(x))=f(\Phi(x))_{+}-f(\Phi(x))\_$ where $f(\Phi(x))_{+}\geq 0$ and $f(\Phi(x))\_\geq 0$ have orthogonal supports, is now given by 

 \[f(\Phi(x))_{+}=\begin{bmatrix}
p_1 f(\Phi(x)) p_1 & 0 & 0 & 0\\
0 & 0 & 0 & 0\\
0 & 0 & 0 & 0\\
0 & 0 & 0 & p_4 f(\Phi(x)) p_4
\end{bmatrix},\]

and  \[f(\Phi(x))\_=\begin{bmatrix}
0 & 0 & 0 & 0\\
0 & -p_2 f(\Phi(x)) p_2 & 0 & 0\\
0 & 0 & -p_3 f(\Phi(x)) p_3 & 0\\
0 & 0 & 0 & 0
\end{bmatrix}.\]

However, the Jordan decomposition of $\Phi(f(x))$ is difficult to describe, but $E(\Phi(f(x)))$ is given by

\[(E(\Phi(f(x))))_{+}=\begin{bmatrix}
p_1 \Phi(f(x)) p_1 & 0 & 0 & 0\\
0 & (p_2 \Phi(f(x)) p_2)_{+} & 0 & 0\\
0 & 0 & (p_3 \Phi(f(x)) p_3)_{+} & 0\\
0 & 0 & 0 & p_4 \Phi(f(x)) p_4
\end{bmatrix},\]

\[(E(\Phi(f(x))))\_=\begin{bmatrix}
0 & 0 & 0 & 0\\
0 & (p_2 \Phi(f(x)) p_2)\_ & 0 & 0\\
0 & 0 & (p_3 \Phi(f(x)) p_3)\_ & 0\\
0 & 0 & 0 & 0
\end{bmatrix}.\]

Now using Lemma \ref{key lemma}
 we can compare the diagonal blocks of $f(\Phi(x))_{\pm}$ and $(E(\Phi(f(x))))_{\pm}$:
 \[p_i \Phi(f(x))p_i\gtrsim p_i f(\Phi(x))p_i, \  \text{for} \ i=1, 4\]

  \[-p_i f(\Phi(x))p_i\gtrsim  (p_i \Phi(f(x))p_i)\_ \  \  \text{for} \ i=2, 3.\]

  Taking direct sum, we get 

  \begin{align*}
  (E(\Phi(f(x)))_{+}&\geq \begin{bmatrix}
p_1 \Phi(f(x)) p_1 & 0 & 0 & 0\\
0 & 0 & 0 & 0\\
0 & 0 & 0 & 0\\
0 & 0 & 0 & p_4 \Phi(f(x)) p_4
\end{bmatrix}\\
&\gtrsim \begin{bmatrix}
p_1 f(\Phi(x)) p_1 & 0 & 0 & 0\\
0 & 0 & 0 & 0\\
0 & 0 & 0 & 0\\
0 & 0 & 0 & p_4 f(\Phi(x)) p_4
\end{bmatrix}=f(\Phi(x))_{+}.
\end{align*}

Thus we get 
\begin{equation}\label{eq-pre-order 1}
  (E(\Phi(f(x)))_{+} \gtrsim f(\Phi(x))_{+}.
\end{equation}
And similarly one obtains \begin{equation}\label{eq-pre-order2}
f(\Phi(x))\_ \gtrsim (E(\Phi(f(x)))\_.
\end{equation}

Thus, equation \ref{eq-pre-order 1} yields $\tau((E(\Phi(f(x)))_{+})\geq \tau(f(\Phi(x))_{+})$ and equation \ref{eq-pre-order2} yields $\tau(f(\Phi(x))\_)\geq \tau((E(\Phi(f(x)))\_).$ Now assume that both the quantities $\tau(\Phi(f(x)))$ and $\tau(f(\Phi(x)))$ are well defined. Using the fact that $E$ is $\tau$-preserving, we get  that $\tau(E(\Phi(f(x))))$ is well defined. Indeed, from the minimality of the Jordan decomposition (see \cite[Section 2.4]{harada-kosaki}) we see that 
 $$\tau(E(\Phi(f(x)))_{\pm})\leq \tau(E(\Phi(f(x))_\pm))=\tau (\Phi(f(x))_{\pm}).$$
Hence we obtain 
\begin{align*}
\tau(E(\Phi(f(x)))&=\tau((E(\Phi(f(x)))_{+})-\tau((E(\Phi(f(x)))\_)\\
&\geq \tau(f(\Phi(x))_{+})- \tau(f(\Phi(x))\_)\\
&=\tau(f(\Phi(x))).
\end{align*}
Now note that $$E(\Phi(f(x)))_+-E(\Phi(f(x)))_-=E(\Phi(f(x)))=E(\Phi(f(x))_+)-E(\Phi(f(x))\_)$$ and hence
$$\tau(E(\Phi(f(x)))_+)+\tau (E(\Phi(f(x))\_))=\tau(E(\Phi(f(x))_+))+\tau(E(\Phi(f(x)))\_),$$
and we get  using the trace preservation of $E$
\begin{align*}\tau(\Phi(f(x)))&=\tau(\Phi(f(x))_+)-\tau(\Phi(f(x))\_)\\&=\tau(E(\Phi(f(x))_+))-\tau(E(\Phi(f(x))\_))\\&=\tau((E(\Phi(f(x)))_{+})-\tau((E(\Phi(f(x)))\_)\\&=\tau(E(\Phi(f(x)))\geq \tau(f(\Phi(x))).
\end{align*}

\section{The proof of the main theorem}

Let $(\mathcal{M}_1, \tau_1)$, $ (\M_2,\tau_2)$ be two tracial von Neumann algebras. For $\omega\in (\M_2)_*\simeq L^1(\M_2,\tau_2)$ let $R_\omega: \M_1\bar\otimes \M_2\to \M_1$ be the right slice operator defined on elementary tensors as $R_\omega(a\otimes b)=a\omega(b)$. Similarly, we define a left slice map $L_w: \M_1\bar\otimes \M_2\to\M_2$, $w\in (\M_1)_*$. The maps are completely bounded normal maps such that $\|L_w\|_{\rm cb}=\|w\|$, $\|R_\omega\|_{cb}=\|\omega\|$. Moreover, if $\omega$ is a positive linear functional (a state) then $R_\omega$ is normal completely positive (and unital).


The same holds for the left slice map.  Then for any $X\in\M_1\bar\otimes\M_2$, $a$, $b\in\M_1$, we have $R_\omega((a\otimes 1)X(b\otimes 1))=aR_\omega(X)b$ (checked on elementary tensor products). 
Assuming $\tau_1(a^*a)<\infty$, i.e. $a\in L^2(\M_1,\tau_1)$, we obtain $a^*R_\omega(X)a$  is in $L^1(\M_1,\tau_1)$, with $$|\tau_1(a^*R_\omega(X))a|\leq \|a\|_2\|R_\omega(X)\|\leq \|a\|_2\|\omega\|\|X\|.$$ Therefore we can define an element $(\tau_1\otimes{\rm id})((a^*\otimes 1)X(a\otimes 1))$ in $\M_2$ by letting
$$\langle (\tau_1\otimes{\rm id})((a^*\otimes 1)X(a\otimes 1)),\omega\rangle=\tau_1(a^*R_\omega(X)a).$$
We have also $\|(\tau_1\otimes{\rm id})((a^*\otimes 1)X(a\otimes 1))\|\leq \|a\|_2\|X\|.$ Note that letting $\omega_{\tau_1,a}(Y)=\tau_1(a^*Ya)$, $Y\in\M_1$, which is a bounded linear functional on $\M_1$, we have $(\tau_1\otimes{\rm id})((a^*\otimes 1)X(a\otimes 1))=L_{\omega_{\tau_1,a}}(X)$. 

\begin{remark}\label{remark}
For $X\in\M_1\bar\otimes\M_2$, we say that $({\rm id}\otimes \tau_2)(X)$ is {\bf defined} whenever $L_\omega(X)\in L^1(\M_2,\tau_2)$ for any $\omega\in(\M_1)_*$ and let $\langle({\rm id}\otimes \tau_2)(X),\omega\rangle:=\tau_2(L_\omega(X))$. By the Uniform Boundedness Principle, $({\rm id}\otimes \tau_2)(X)$ is indeed a bounded linear functional on $(\M_1)_*$ and hence is an element of $\M_1$.

In Theorem \ref{main theorem}, by (\ref{eq-main}) being defined we  mean that $({\rm id}\otimes \tau_2)(f(H))$ is defined and show that in this case the left hand side is defined in the sense of  Definition \ref{def-Jordan}. 
\end{remark}
\begin{lemma}\label{tau2}
Let $a\in L^2(\M_1,\tau_1)$ and  $X\in\M_1\bar\otimes\M_2$. If  $({\rm id}\otimes \tau_2)(X)$ is defined, then $(\tau_1\otimes{\rm id})((a^*\otimes 1)X(a\otimes 1))\in L^1(\M_2,\tau_2)$ and $$\tau_2((\tau_1\otimes{\rm id})((a^*\otimes 1)X(a\otimes 1)))=\tau_1(a^*({\rm id}\otimes\tau_2)(X)a).$$
\end{lemma}
\begin{proof}
    Consider $\omega\in (\M_1)_*$, given by $\omega(M)=\tau_1(a^*Ma)$, $M\in \M_1$. Then for any $w\in (\M_2)_*$, as
    \begin{eqnarray*} 
    &&\langle L_\omega(X),w\rangle=\langle X,\omega\otimes w\rangle=\omega(R_w(X))\\&&=\tau_1(a^*R_w(X)a)=\langle (\tau_1\otimes{\rm id})(a^*\otimes 1)X(a\otimes 1)),w\rangle,
    \end{eqnarray*}
    we obtain 
    $L_\omega(X)=(\tau_1\otimes{\rm id})((a^*\otimes 1)X(a\otimes 1))\in L^1(\M_2,\tau_2)$. 
    
    On the other hand, $\langle ({\rm id}\otimes \tau_2)(X),\omega\rangle =\tau_1(a^*({\rm id}\otimes \tau_2)(X)a)$, giving
    $$\tau_2(L_\omega(X))=\tau_2((\tau_1\otimes{\rm id})(a^*\otimes 1)X(a\otimes 1))=\tau_1(a^*({\rm id}\otimes \tau_2)(X)a).$$
\end{proof}

\textbf{Proof of Theorem \ref{main theorem}:}

Define a map $\Phi: \M_1\bar{\otimes} \M_2\rightarrow \M_2$ by 
    \[\Phi(X)=(\tau_1\otimes {\rm id}) ((a^*\otimes 1)X (a\otimes 1)).\]
Note that the condition $\tau_1(a^*a)= 1$ makes the map well-defined, completely positive and unital as the slice map is $L_{\omega_{\tau_1,a}}$ with $\omega_{\tau_1,a}(y)=\tau_1(a^*ya)$, $y\in\M_1$; see the above arguments.
    
    Using Theorem \ref{improved Petz thm} we get for any self-adjoint $H$
    \[\tau_2 (f (\Phi(H)))\leq \tau_2 (\Phi(f(H))),\]
    whenever both sides are defined. 
    
    As  $\Phi(f(H))\in L^1(\M_2,\tau_2)$,  the right hand side is defined and finite. Moreover, if $E$ is as in the proof of Theorem  \ref{improved Petz thm}, $E(\Phi(f(H)))\in L^1(\M_2,\tau_2)$ and $\tau_2 (\Phi(f(H)))=\tau_2 (E(\Phi(f(H))))$ and $\tau_2((E(\Phi(f(H)))_+)$, 
    $\tau_2(E(\Phi(f(H)))_-)$ are finite. It now follows from the proof of Theorem \ref{improved Petz thm} that $\tau_2(f(\Phi(H))_+)\leq\tau_2(E(\Phi(f(H)))_+)<\infty$ and hence the left hand side is defined, but can take the value $-\infty$.  
This yields

\[\tau_2 f (\tau_1\otimes {\rm id})[(a^*\otimes 1)H (a\otimes 1)]\leq \tau_2 (\tau_1\otimes {\rm id})[(a^*\otimes 1)f(H) (a\otimes 1)].\]
Now using Lemma \ref{tau2} to the RHS of the above inequality we get 

\[\tau_2 f (\tau_1\otimes {\rm id})[(a^*\otimes 1)H (a\otimes 1)]\leq  \tau_1 (a^* ({\rm id}\otimes \tau_2)(f(H))a),\]
which yields the desired result.

\section{Jensen's inequality for states (non-tracial case)} 
In this section we show a similar Jensen inequality  with general, not necessarily, tracial states. 
We need the notion of operator convexity which is a stronger notion than  that of ordinary convexity. Note that a continuous function $f$ on an interval $I$ is said to be \textbf{operator convex} if 
$f(\lambda x+ (1-\lambda) y)\leq \lambda f(x)+(1-\lambda) f(y)$ holds whenever $\lambda\in [0, 1]$ and for any self-adjoint operators $x, y$ with spectra included $I$. Operator convexity is stronger than ordinary convexity \cite{Bhatia1997}.

\begin{theorem}\label{main theorm states}
    Let $\mathcal{M}_1, \M_2$ be two von Neumann algebras and two normal states $\rho_1, \rho_2$  acting on $\M_1, \M_2$ respectively. Suppose that $H$ is a self-adjoint element in $\M_1\bar{\otimes} M_2$. Then for any contraction $a\in \M_1$ and any operator convex function $f$ defined on the spectrum of $H$ we have
    \begin{equation}\label{eq-1}
    \rho_2 f [  (\rho_1\otimes {\rm id})((a^*\otimes 1) H (a\otimes 1))]\leq  \rho_1 (a^* ({\rm id}\otimes \rho_2)(f(H))a).
    \end{equation}
\end{theorem}
\begin{proof}
    We note that for states $\rho_1, \rho_2$, the corresponding left and right slice maps are unital completely positive maps, that is, $L_{\rho_1}: \M_1\bar\otimes \M_2\rightarrow \M_2$ defined by $L_{\rho_1}(x\otimes y)=(\rho_1\otimes {\rm id}) (x\otimes y)=\rho_1(x)y$ and similarly $R_{\rho_2}: \M_1\bar\otimes \M_2\rightarrow \M_1$ defined by
    $R_{\rho_2}(x\otimes y)=({\rm id}\otimes \rho_2)(x\otimes y)=x\rho_2(y)$.
    Note that these maps are unital because the $\rho_1, \rho_2$ are states, so we have $\rho_1(1)=1=\rho_2(1)$, where we abused our notations slightly and denoted the identity element of both the algebras as $1$ and the number appears in the middle is the numeric $1$.
    
    Note that for any ucp map $\Phi$ and an operator convex function $f$ on an interval $I$, we have $\Phi(f(H))\geq f(\Phi(H))$ for every self-adjoint $H$ (see  \cite{choi1974, davis} ) such that spectrum of $H$ is inside $I$.
    Now calling $X:=(a^*\otimes 1)  H (a\otimes 1)$ which is a self-adjoint element, we have for the left-hand side of the equation \ref{eq-1}
    \begin{equation}\label{eq-2}
    \rho_2 f [ (\rho_1\otimes {\rm id})((a^*\otimes 1)  H (a\otimes 1))]=\rho_2 ( f (L_{\rho_1}(X)))\leq \rho_2 (L_{\rho_1}(f(X))).
    \end{equation}
    
    Now note that for any $A\in \M_1\bar\otimes M_2$ we have $$\rho_2 (L_{\rho_1} (A))=(\rho_1\otimes \rho_2)(A)=\rho_1(R_{\rho_2}(A)).$$ 

    Now using the above commutation relation we see the right hand side of the equation \ref{eq-2}
\begin{equation}\label{eq-3}
 \rho_2 (L_{\rho_1}(f(X)))=\rho_1 (R_{\rho_2}(f[(a^*\otimes 1)  H (a\otimes 1)]))\leq \rho_1 (R_{\rho_2} ((a^*\otimes 1)f(H)(a\otimes 1))),
\end{equation}    
where we have used Jensen's operator inequality (contractive version) from Hansen and Pederson (see Corollary 2.3 in \cite{Hansen-Pedersen}): $$f((a^*\otimes 1)H (a\otimes 1))\leq (a^*\otimes 1)f(H)(a\otimes 1),$$
and positivity of the map $\rho_1 R_{\rho_2}$.
As for any $a\in \M_1$ and $A\in \M_1\bar{\otimes} \M_2$,
\[R_{\rho_2} ((a^*\otimes 1)A(a\otimes 1))=a^*R_{\rho_2}(A)a,\]
 we get that the right hand side of equation \ref{eq-3} is 
\[\rho_1 R_{\rho_2} ((a^*\otimes 1)f(H)(a\otimes 1)))=\rho_1 [a^*R_{\rho_2}(f(H))a].\]
Putting all together we have 
\[\rho_2 f [ (\rho_1\otimes {\rm id})((a^*\otimes 1)  H (a\otimes 1))] \leq \rho_1 [a^*({\rm id}\otimes \rho_2)(f(H))a].\]

\end{proof}
\begin{remark}
    We remark that normal states  $\rho_1$, $\rho_2$ in Theorem \ref{main theorm states} can be replaced by unital normal completely positive maps to have the following operator inequality:
     $$
    ({\rm id} \otimes\rho_2)( f [  (\rho_1\otimes {\rm id})((a^*\otimes 1) H (a\otimes 1))])\leq  (\rho_1\otimes {\rm id}) (a^*\otimes 1 ({\rm id}\otimes \rho_2)(f(H))a\otimes 1).
    $$
    Working with the scalar case it would be interesting to know whether the statement of Theorem \ref{main theorm states} holds with continuous convex function $f$ instead of operator convex one. 
\end{remark}

\section{Outlook}

The finite-dimensional techniques used in the proof of Theorem 1 can be extended to the infinite-dimensional setting under certain assumptions on the discreteness of the spectrum, as noted in \cite{Carlen-Frank=Larson}. Under these conditions, the authors obtained results on eigenvalue asymptotics. We expect that our theorem will allow one to treat more general situations beyond this setting.

Trace inequalities, and operator Jensen’s inequality in particular, have played an important role in the study of quantum entropy of quantum systems (see \cite{effros, carlen-book}). It would be interesting to investigate whether the results of this paper have implications for problems arising in quantum information theory.
\section{Acknowledgment} We thank Simon Larson for many interesting discussions during the course of this work. We also thank the anonymous referee for many helpful suggestions that improved the exposition of the article. M.R acknowledges support from
the Wallenberg Center for Quantum Technology
(WACQT) and L.T was supported by the Swedish Research Council project grant 2023-04555.

\bibliography{references}
\bibliographystyle{alpha}
\end{document}